\documentclass[12pt, reqno]{amsproc}
\usepackage{graphicx, enumerate}
\usepackage[margin=1in]{geometry}
\usepackage{amsmath}
\usepackage{amsfonts}
\usepackage{amssymb}
\usepackage{amsfonts,amsmath,amssymb}
\usepackage{amsthm}
\newtheorem{lemma}{Lemma}
\newtheorem{theorem}{Theorem}

\newtheorem{defi}{Definition}
\newtheorem{eg}{Example}
\title[Order of growth of solutions of linear complex differential equations when the coefficients have same order]{Order of growth of linear complex differential equations when the coefficients have same order}
\author[Naveen Mehra and S. K. Chanyal]{Naveen Mehra and S. K. Chanyal}
\address{Naveen Mehra, Department of Mathematics, Kumaun University, D.S.B. Campus, Nainital-263001, Uttarakhand, India}
\email{naveenmehra00@gmail.com}
\address{S. K. Chanyal, Department of Mathematics, Kumaun University, D.S.B. Campus, Nainital-263001, Uttarakhand, India}
\email{skchanyal.math@gmail.com}
\thanks{The research work of the first author is supported by research fellowship from Council of Scientific and Industrial Research (CSIR), New Delhi.} 
\begin{document}
\maketitle 
\begin{abstract}
This article contains the theorems which shows that when $A(z)=h_1(z)e^{P_1(z)}$ and $B(z)=h_0(z)e^{P_0(z)}$ are of same order,then all the non-trivial solutions of equation $f''+A(z)f'+B(z)f=0$ are of infinite order. Moreover we extend these results to higher order differential equations.
\end{abstract}
\section{Introduction and statement of main results}
Consider the differential equation, 
\begin{equation}\label{cde}
f''+A(z)f'+B(z)f=0
\end{equation}
where $A(z)$ and $B(z)(\not\equiv0)$ are entire functions. Then, it is well known that all solutions of equation \eqref{cde} are entire functions(see \cite{herold}).
\begin{defi}
The order of growth $\rho(f)$ and lower order of growth $\mu(f)$ of an entire function $f$ are defined by\begin{center}
 $\displaystyle{\rho(f)=\limsup\limits_{r \to \infty}\displaystyle{\frac{\log^+\log^+{M(r,f)}}{\log{r}}}}$
 $\displaystyle{\mu(f)=\liminf\limits_{r \to \infty}\displaystyle{\frac{\log^+\log^+{M(r,f)}}{\log{r}}}}$ \end{center}
 respectively, where $M(r,f)= max_{|z|=r}|f(z)|$.
\end{defi} 
\begin{defi}
Exponent of convergence of zeros of function $f$ is defined by
\begin{center}
 $\displaystyle{\lambda(f)= inf\lbrace\alpha>0\displaystyle{|}\sum\limits_{n=1}^{\infty}|a_n|^{-\alpha}<\infty\rbrace},$
 \end{center}
where $(a_n)_{n\in\mathbb{N}}$ is a sequence of non-zero complex numbers(not-necessarily distinct) with $\lim_{n\to\infty}a_n=\infty$.\end{defi}
 If $A(z)$ and $B(z)$ are polynomials, then all solutions of  equation \eqref{cde} are of finite order(see \cite{wittich}). Moreover, if either of $A(z)$ and $B(z)$ or both are entire transcendental function, then equation \eqref{cde} have atleast one solution of infinite order(see \cite{frei}). Following classical result is summary of \cite{Gundersen} and  \cite{heller}.
\begin{theorem}
Suppose that $A(z)$ and $B(z)$ are entire functions satisfying one of the following conditions:
\begin{enumerate}
\item[(i)]\cite{Gundersen} $\rho(A)<\rho(B)$; 
\item[(ii)]\cite{Gundersen} $A(z)$ is a polynomial and $B(z)$ is a transcendental entire function;
\item[(iii)]\cite{heller} $\rho(B)<\rho(A)\leq\frac{1}{2}$.
\end{enumerate}
Then every non-trivial solution of equation (1) are of infinite order.
\end{theorem}
Following are the examples of differential equations having finite order solutions when $\rho(B)\leq \rho(A)$.
\begin{eg}
In the differential equation $f''+(e^{z}-1)f'+e^zf=0$, $f(z)=e^z$ is a non-trivial solution of order $1$ and order of both coefficients are equal.
\end{eg}
\begin{eg}
The differential equation $f''+(ze^z+z+1)f'+zf=0$ has a solution $f=e^{-z}+1$.
\end{eg}
Above examples shows the differential equation \eqref{cde} has finite order solution when $\rho(B)\leq\rho(A)$. 
Consider $A(z)=h_1(z)e^{P_1(z)}$ and $B(z)=h_0(z)e^{P_0(z)}$, where $h_1(z)$ and $h_0(z)$ are transcendental entire functions; $P_1(z)$ and $P_0(z)$ are polynomials satisfying $\rho(h_1)<\rho(h_0)$. Gundersen \cite{gundersen2} proved that all non-trivial solutions are of infinite order when $\deg P_1\neq \deg P_0$. What if $\deg P_1 = \deg P_0$? 
\begin{eg}
$f''(z)-z^2e^z+ze^zf=0$ has finite order non-trivial solution $f(z)=2z$.
\end{eg}
In example 3, we have $\rho(A)=\rho(B)=\deg(z)=1$, where $\lambda(A)<\rho(A)$ and $\lambda(B)<\rho(B)$. Let 
\begin{equation}\label{eqpolyP}
P_1(z)= a_nz^n + a_{n-1} z^{n-1}+\cdots + a_0
\end{equation}
 and 
\begin{equation}\label{eqpolyQ}
P_0(z)= b_nz^n + b_{n-1} z^{n-1}+\cdots + b_0
\end{equation} 
where $a_n\neq 0$ and $b_n\neq 0$. Kwon\cite{kwon} considered  this case with some restrictions and gives condition when order of solution is infinite. He proved the following:

\begin{theorem}\cite{kwon} 
Consider the equation
\begin{equation}
\label{cde2}
f''+h_1(z)e^{P_1(z)}f'+h_0(z)e^{P_0(z)}f=0,
\end{equation}
where $P_1(z)$ and $P_0(z)$ are non-constant polynomials as given in \eqref{eqpolyP} and \eqref{eqpolyQ} with $a_n\neq 0$, $b_n\neq 0$ and $h_1(z)$, $h_2(z)\neq 0$ are entire functions satisfying $\rho(h_1)< degP_1$ , $\rho(h_0)<degP_0$. Then all non-trivial solutions are of infinite order if any of the following holds:
\begin{enumerate}[(i)]
\item $\arg a_n\neq \arg b_n$ or $a_n=cb_n$ with $0<c<1$,
\item $a_n=b_n$ and $\deg (P_1-P_0)=m\geq 1$, $\rho(h_1)<m$ and $h_0(z)<m$,
\item $a_n-cb_n$ with $c>1$ and $\deg(P_1-cP_0)=m\geq 1$, $\rho(h_1)<m$
and $0<\rho(h_0)<\frac{1}{2}$ or
\item $a_n=cb_n$ with $c\geq 1$, $P_1(z)-cP_0(z)$ be a constant and
$\rho(h_1)<\rho(h_0)<\frac{1}{2}$.
\end{enumerate}
\end{theorem}
In this article, we shall like to expand Kwon's result when $P_1(z)-cP_0(z)=constant$ and $c\geq 1$. His result has a limitation on the $\rho(h_0)<\frac{1}{2}$. In our first result we consider $\rho(h_1)<\mu(h_0)<\frac{1}{2}$. The statement of our first result is following:
\begin{theorem}\label{mainthm1}
Let $P_1(z)$ and $P_0(z)$ be non-constant polynomials as given in \eqref{eqpolyP} and \eqref{eqpolyQ} with $a_n\neq 0$, $b_n \neq 0$ and $h_1(z)$, $h_0(z)(\not\equiv 0)$ be an entire functions satisfying $\rho(h_1)< degP_1$ and $\rho(h_0) < degP_0$. Let  $P_1(z) - cP_0(z)$ be a constant  with $c\geq 1$ and $\rho(h_1)<\mu(h_0)\leq\frac{1}{2}$. Then all non-trivial solutions of \eqref{cde2} are of infinite order.
\end{theorem}
In Theorem \ref{mainthm1}, we restrict $\mu(h_0)$ to be less than half. If we consider $h_1(z)$ to have Fabry gaps or multiply connected Fatou component, we consider $\mu(h_0)$ cannot be less than half. 
\begin{defi}
An entire function $f=\sum_{n=0}^{\infty}a_{\lambda_n}z^{\lambda_n}$ is said to have Fabry gaps if $\frac{\lambda_n}{n}\to\infty$.
\end{defi}
There are many research articles in which authors used coefficient having a Fabry gaps. Long\cite{long}, Kumar and Saini\cite{kumarsaini} also used Fabrys gap in their articles. 

\begin{theorem}\label{mainthm2}
Let $P_1(z)$ and $P_0(z)$ be non-constant polynomials as given in \eqref{eqpolyP} and \eqref{eqpolyQ} with $a_n\neq 0$, $b_n \neq 0$ and $h_1(z)$, $h_0(z)(\not\equiv 0)$ be an entire functions satisfying $\rho(h_1)< degP_1$ and $\rho(h_0) < degP_0$. Let  $P_1(z) - cP_0(z)$ be a constant  with $c\geq 1$ and let $h_1(z)$ has Fabry gaps with $\mu(h_0)>\rho(h_1)$. Then every non-constant solution $f$ of \eqref{cde2} has infinite order.
\end{theorem}
\begin{defi}
For the entire function, the set of family of iterates $\{f^n:n\in\mathbb{N}\}$ is well defined and forms a normal family; is called Fatou set.  
\end{defi}
Interested reader may visit \cite{devaney} for detailed analysis about Fatou set. The complement of Fatou set is known as Julia set. G. Zhang and L. Yang\cite{zhangyang} considered entire function having multiply connected Fatou component in his paper.
\begin{theorem}\label{mainthm3}
Let $P_1(z)$ and $P_0(z)$ be non-constant polynomials as given in \eqref{eqpolyP} and \eqref{eqpolyQ} with $a_n\neq 0$, $b_n \neq 0$ and $h_1(z)$, $h_0(z)(\not\equiv 0)$ be an entire functions satisfying $\rho(h_1)< degP_1$ and $\rho(h_0) < degP_0$. Let  $P_1(z) - cP_0(z)$ be a constant  with $c\geq 1$ and $h_1(z)$ has multiply connected Fatou component with $\rho(h_1)\leq\mu(h_0)$. Then all non-trivial solutions of \eqref{cde2} are of infinite order.
\end{theorem}
\section{Preliminary Lemma} 
This section is devoted to preliminaries results which are used to prove our main theorems when a transcendental entire function $f$ is of finite order Gundersen\cite{Gundersen} gave a well known result. This result is used extensively.
\begin{lemma}\cite{Gundersen}\label{Gundersen} Let $f$ be a transcendental entire function with $\rho<\infty$, let $\Gamma = \lbrace(k_1,j_1), \\(k_2, j_2),\ldots,(k_n, j_n)\rbrace$ denotes finite set of distinct pairs of integers that satisfy $k_i > j_i \geq 0$, for $i = 1, 2,\ldots ,n$, and let
$\epsilon > 0$ be a given constant. Then the following three statements holds:\begin{enumerate}
\item[(i)] there exists a set $E_1 \subset [0, 2\pi)$ that has linear measure $0$, such that if $\psi_0 \in [0, 2\pi) -E_1$, then there is a constant $R_0 = R_0(\psi_0) > 0$ so that for all z satisfying $\arg z = \psi_0$ and $|z| \geq R_0$ and for all $(k, j) \in \Gamma$, we have
\begin{equation}\label{f"byf}
\left|\frac{f^{(k)}(z)}{f^{(j)}(z)}\right| \leq |z|^{(k-j)(\rho-1+\epsilon))}.
\end{equation}
\item[(ii)] there exists a set $E_2 \subset (1,\infty)$ that has finite logarithmic measure, such that for all z with $|z| \notin E_2 \cup[0, 1]$ and for all $(k, j) \in \Gamma$,
the inequality \eqref{f"byf} holds.
\item[(iii)] there exists a set $E_3 \subset [0,\infty)$ that has finite linear measure, such that for all z with $|z| \notin E_3$ and for all $(k, j) \in \Gamma$, we have 
\begin{equation}
\left|\frac{f^{(k)}(z)}{f^{(j)}(z)}\right| \leq |z|^{(k-j)(\rho+\epsilon)}.
\end{equation}
\end{enumerate}
\end{lemma}
Besicovitch's $\cos\pi\rho$ theorem provides the condition for a function $f$ of order less than 1. It says that $\log  m(r,f)>(cos \pi\rho')\log M(r,f)$ when  $\rho<\rho' <1$, $|z|=r$ belongs to a set of upper density greater than $1-\frac{\rho}{\rho'}$, where $m(r,f)=\min\{|f|: |z|=r\}$ and $M(r,f)=\max\{|f|: |z|=r\}$. Barry\cite{barry} proved this result for $\mu <1$.
\begin{lemma}\cite{barry}\label{lcospirho}
Assume that $f$ is an entire function with $0\leq \mu(f) < 1$. Then, for every
$\beta\in (\mu(f), 1)$, $m(r,f) > M(r,f)\cos\pi\alpha$ for $|z|=r\in F\subset [1,\infty)$, where $\overline{\log dens}F \geq 1-\frac{\mu(f)}{\beta}$.
\end{lemma}
Combining Theorem 1 from \cite{Fuchs} and Lemma 2.2 from \cite{Ishizaki}, we get the following lemma. It is proved by Long\cite{long}. It is the property of finite order entire functions having Fabry gaps.
\begin{lemma}\cite{long}\label{fabry}
Assume $f(z) =\sum\limits_{n=0}^\infty a_{\lambda_n}z^{\lambda_n}$ be an entire function of finite order of growth with Fabry gaps, and let $u(z)$ be an entire function with $\rho(u)  \in (0,\infty)$. Then for any given $\epsilon \in (0, \varsigma)$, where $\varsigma=\min(1,\rho(u))$, there exists a set $G \subset (1,\infty)$ satisfying $\overline{\log dense}G \geq \xi$, where $\xi \in (0, 1)$ is a constant such that for all $|z| = r \in G$, $$\log M(r, u) > r^{\rho(u) -\epsilon},\   \log m(r, g) > (1-\xi) \log M(r, g),$$ where $M(r, u) = \max\lbrace|u(z)| : |z| = r\rbrace ,\ m(r, g) = \min\lbrace|g(z)| : |z| = r\rbrace$ and $M(r, g) = \max\{ |g(z)| : |z| = r\}$.
\end{lemma}
Every multiply connected Fatou component of a transcendental meromorphic function must be Baker wandering domain(proved by Baker\cite{baker}) and if a transcendental meromorphic function has Baker wandering domain, then $\mathcal{J}$ has only bounded components(see \cite{zheng} Remark on page $25$).  So, following lemma was proved by Zheng\cite{zheng}(see Corollary 1), also holds for transcendental meromorphic function having multiply connected Fatou component. 
\begin{lemma}\cite{zheng}\label{zheng}
Assume that f is a transcendental meromorphic function with at most finitely many poles. If $\mathcal{J}(f)$ has only bounded components, then for any complex number $z$,  there exist a constant $0<\gamma<1$ and two sequences $\{r_n\}$ and $\{R_n\}$ of positive numbers with $r_n\to\infty$ and $R_n/r_n\to\infty(n\to\infty)$ such that
$$M(r,f)^\gamma\leq L(r,f), r\in H,$$ where $H=\cup_{n=1}^\infty\{r:r_n<r<R_n\}.$  
\end{lemma}
\section{Proof of the main theorems}
\subsection*{Proof of Theorem \ref{mainthm1}.}
\begin{proof}
Let $\rho(f)<\infty$. Then by Lemma \ref{Gundersen}, there exist a set $E\subset [0,\infty)$ having finite linear measure such that for $r\not\in E$ 
\begin{equation}\label{eqffiniteorder}
\left|\frac{f^{(k)}(re^{\iota\theta})}{f(re^{\iota\theta})}\right| \leq r^{k\rho(f)}.
\end{equation}
Since $\mu(h_0)<1/2$, by Lemma \ref{lcospirho} for $r\in F$, where $\overline{\log dens}F \geq 1-\frac{\mu(f)}{\beta}$, we have 
\begin{equation}\label{eqcospirhothm3}
M(r,h_0)\cos\pi\alpha\ < |h_0(re^{\iota\theta})|. 
\end{equation}
By definition of order, we have for $\epsilon>0$ 
\begin{equation}\label{orderh_1}
|h_1(re^{r\iota\theta})|\leq \exp{r^{\rho(h_1) + \epsilon}},
\end{equation}
holds for all $r\geq r_0$, where $r_0>0$.
From \eqref{cde},we have
\begin{equation}\label{eqh2incdethm3}
|h_0(re^{\iota\theta})e^{(1-c)P_0(re^{\iota\theta})}|\leq |e^{-cP_0(re^{\iota\theta})}|\left|\frac{f''(re^{\iota\theta})}{f(re^{\iota\theta})}\right|+|h_1(re^{\iota\theta})e^{P_1(re^{\iota\theta})-cP_0(re^{\iota\theta})}|\left|\frac{f'(re^{\iota\theta})}{f(re^{\iota\theta})}\right|.
\end{equation}
Using \eqref{eqffiniteorder}, \eqref{eqcospirhothm3}  and \eqref{orderh_1} in (\ref{eqh2incdethm3}), for very large  $r\in F\setminus E$ on a curve $C$ be a tending to $\infty$ such that $ReP_0(z)=0$, we get
$$
M(r,h_0)\cos\pi\beta < (1+ \exp{r^{\rho(h_1) + \epsilon}})r^{2\rho(f)},
$$
which is a contradiction, since $\rho(h_1)<\mu(h_0)$.
\end{proof}
\subsection*{Proof of Theorem \ref{mainthm2}.}
\begin{proof}
Assume that solution $f$ of equation \eqref{cde} to be of finite order. Using  Lemma \ref{fabry} on transcendental entire function $h_0(z)$ having Fabry gaps,  for $r\in G$, where $\overline{\log dens}G \geq \xi$, for $\xi\in(0,1)$ we have 
\begin{equation}\label{eqfabry}
M(r,h_0)^{1-\xi} < |h_0(re^{\iota\theta})|. 
\end{equation}
From \eqref{cde}, we have
\begin{equation}\label{eqh2incdethm4}
|h_0(re^{\iota\theta})e^{(1-c)P_0(re^{\iota\theta})}|\leq |e^{-cP_0(re^{\iota\theta})}|\left|\frac{f''(re^{\iota\theta})}{f(re^{\iota\theta})}\right|+|h_1(re^{\iota\theta})e^{P(re^{\iota\theta})-cP_0(re^{\iota\theta})}|\left|\frac{f'(re^{\iota\theta})}{f(re^{\iota\theta})}\right|.
\end{equation}
Using \eqref{eqffiniteorder}, \eqref{orderh_1}    and \eqref{eqfabry} in (\ref{eqh2incdethm4}), for very large  $r\in G\setminus E$ on a curve $C$ be a tending to $\infty$ such that $ReP_0(z)=0$, we get
$$M(r,h_0)^{1-\xi} < (1+ \exp{r^{\rho(h_1) + \epsilon}})r^{2\rho(f)}
,$$
which is a contradiction.
\end{proof}
\subsection*{Proof of Theorem \ref{mainthm3}.}
\begin{proof}
Suppose that $\rho(f)<\infty$. Considering $h_0$ to be an entire function having multiply connected Fatou component, by Lemma \ref{zheng} for $r\in H$, where $H=\cup_{n=1}^\infty\{r:r_n<r<R_n\}$ , we have 
\begin{equation}\label{eqzheng}
M(r,h_0)^\gamma < |h_0(re^{\iota\theta})| 
\end{equation}
From \eqref{cde}, we have
\begin{equation}\label{eqh2incdethm5}
|h_0(re^{\iota\theta})e^{(1-c)P_0(re^{\iota\theta})}|\leq |e^{-cP_0(re^{\iota\theta})}|\left|\frac{f''(re^{\iota\theta})}{f(re^{\iota\theta})}\right|+|h_1(re^{\iota\theta})e^{P_1(re^{\iota\theta})-cP_0(re^{\iota\theta})}|\left|\frac{f'(re^{\iota\theta})}{f(re^{\iota\theta})}\right|.
\end{equation}
Using \eqref{eqffiniteorder}, \eqref{orderh_1}  and \eqref{eqzheng}  in (\ref{eqh2incdethm5}), for very large  $r\in H\setminus E$ on a curve $C$ be a tending to $\infty$ such that $ReP_0(z)=0$, we get
$$
M(r,h_0)^\gamma < (1+ \exp{r^{\rho(h_1) + \epsilon}})r^{2\rho(f)}
,$$
which is not possible.
\end{proof}

\end{document}